\documentclass[11pt]{amsart}

\usepackage{amsmath}
\usepackage{hyperref}
\usepackage{xy}
\usepackage{amscd}
\newfont{\sheaf}{eusm10 scaled\magstep1}

\DeclareMathOperator{\Aut}{Aut}

\begin{document}
$$ $$ $$ $$ $$ $$ \bigskip
\title{On the life and scientific work of Gino Fano}

\thanks{ }
\bibliography{mybib}{}
\bibliographystyle{LaTex}
\author[A. Collino]{ Alberto Collino}
\address{A. Collino Universit\'a di Torino, Dipartimento di Matematica, via Carlo Alberto 10 
\indent 10123 Torino, Italy}
 \email{{\tt
alberto.collino@unito.it}}
\author[A. Conte]{Alberto Conte}
\address{ A. Conte Universit\'a di Torino Dipartimento di Matematica, via Carlo Alberto 10 
\indent 10123 Torino, Italy}
 \email{{\tt
alberto.conte@unito.it}}

\author[A. Verra]{Alessandro Verra}
\address{A. Verra Universit\'a Roma Tre, Dipartimento di Matematica e Fisica, Largo San Leonardo Murialdo 1 
\indent 00146 Roma, Italy}
 \email{{\tt
verra@mat.unirom3.it}}
\maketitle
\section  {Gino Fano: a necessarily brief presentation}
 It is evident to most mathematicians, though any estimate seems likely impossible, that the name of Gino Fano is among the most cited ones, in contemporary
geometry and mathematics, nowadays as well as yesterday. \par Furthermore, various fortunate geometric objects bear today his name. That's the case for Fano plane in combinatorics and projective geometry, for Fano varieties in algebraic geometry or for the Fano surface of lines of a cubic hypersurface of a 4-dimensional projective space. \par 
Names, or the number of citations, are however irrelevant with respect to the personality we are dealing with. They play, as it always must be, a secondary role to fully understand a scientific life, the importance of his effective realizations, his heritage. \par
 Fano, during all his life, was a prominent mathematician of his time and a major protagonist of the Italian
School of Geometry. He played a crucial role in the evolution of several fields of Geometry and in strengthening the relations between Italy and the most relevant mathematical and scientific centers of Europe. \par
Fano's life developed between 1871 and 1952.  From the last decade of XIX century, to the period immediately after second World War, his achievements are concerned with the most important topics in Geometry. 
A very synthetic vision of all this is what we offer in the next introductory pages. \par On the other hand, it cannot be the purpose of the present foreword to analyze, with the necessary ampleness and details, Fano's scientific life, nor his influence in today's Geometry and his vaste scientific heritage. All these aspects cannot be confined, specially in the case we are considering, to a small amount of pages of a short, or longer, introduction. \par
The only analysis of Fano's contributions to Algebraic Geometry, even restricting it to one of his favorite topics, like for instance the projective classification of all families of Fano threefolds, is far beyond the limits of an introduction. \par Even more, the same can be said about a historical analysis of Fano's points of view, interests and contributions in a very wide geometric area, which does not exactly coincide with birational classification of algebraic varieties nor with current algebraic geometry. \par Concerning this, the realization of a complete scientific biography of Gino Fano, standing in the very rich field of existing contributions to the study of Italian School of Geometry, appears as a mandatory task whose interest and importance are maximal.  \par
Here our goals are limited\footnote{\tiny \bf As a consequence: no bibliography of papers dedicated to Fano's biography or work is present in this introduction. The  citations in the footnotes only reflect immediate personal choices, without a systematic attempt of completeness. We apologize for every omission.} as follows: in the concluding part of this section we outline, as far as possible, the principal directions and the arguments of Fano's scientifc work or, at least, some jewels from his work. \par The last two sections  are timelines.   The first one reports the main events from Fano's life and scientific activity. Some references to main, or mainly related, historical events are included. \par The second one is a chronologically ordered list of Fano's publications. Every item from this list is quoted in the text by its year, followed by a number denoting its position in the year's list. Starting from 1890, the year of the first publication due to him\footnote{\tiny The translation of Klein's Erlangen Program: see \cite{1}},  the last reprinted book appeared in 1958\footnote{\tiny \it Lezioni di geometria analitica e proiettiva, \rm a famous textbook: see \cite{124}}.  
 \bigskip  \par We hope that the contents of these essential chronologies could offer basic informations, or precisions, to any person interested to the work of Gino Fano and his life.   However two specially nice and vivid descriptions of Gino Fano, as a man and as a scientist, follow this introduction: a memory written by his son Robert Fano and the obituary written by his friend and colleague Alessandro Terracini\footnote{\tiny A. Terracini \it Necrologio di Gino Fano \rm Boll. U.M.I., 7 (1952) 485-490}   \bigskip \par
Like for the greatest painters, the works of great mathematicians often configure several  periods and themes, according to the variations of the intellectual life of their authors and to historical conditions. The works and the life of Gino Fano are a good example of this general principle. \par Two cultural frames strongly contributed to his mathematical education and strongly influenced his scientific personality. To make clear what are these frames, it will be enough to mention two academic cities, where the study of Geometry was flourishing in the last decades of XIX century: \par the city of Torino, where the golden age of Italian Algebraic Geometry was approaching its best moments, and the city of Goettingen where Felix Klein was living since 1886\footnote{\tiny About this  interplay cfr. \it From Turin to Goettingen dialogues and correspondences \rm by E. Luciano and C.S. Roero Boll. storia sci. mat. 32 (2012) 1-232}. \par  Fano, born in the same day and year of Federigo Enriques (January 5 1871) is definitely associated to this golden age since his very beginning. He came to the University of Torino as a student in 1888, obtaining the degree in Mathematics in 1892 under the direction of Corrado Segre. \par He became immediately part of the group of algebraic geometers working in Torino at that time and he started in Torino his academic career, being assistant of Enrico D' Ovidio. It is also due to those years the beginning of a very deep scientific interchange, as well as a permanent friendship, with Guido Castelnuovo. \par
Remarkably, a direct contact with further and very strong perspectives on geometrical studies was made possible to Fano. The reason of this is his visit in the University of Goettingen under the invitation of Felix Klein. Fano spent there about one year, in the period 1893-94. Fano, as Alessandro Terracini writes in the above mentioned memory, had already published  in 1890 a translation  from German to Italian of Klein's Erlangen Program\footnote{\tiny F. Klein \it Vergleichende Betrachtungen Ÿber neuere geometrische Forschungen \rm Verlag von A. Deichert, Erlangen, 1872}, following a suggestion of Corrado Segre. \par That was the first translation from German of this fundamental paper of Klein.  Thus the familiarity with the ideas of Klein, a very fruitful personal interchange with him, as well as Goettingen mathematical atmosphere, gave to the young Fano a further cornerstone, to be connected to the former cornerstone made by his `Italian geometric features', so to build the strong and original personality he was in all the geometric fields, during all his life. Building on such cornerstones, themes and periods of Fano's scientific life evolved in several directions. \par The influence of Felix Klein, as well as  Fano's relations to German School, are probably  reflected by the very high number of his contributions where the general notion of group of geometric transformations takes a central place. This happens in a wide multiplicity of situations. \par His expository memory on continuous groups and geometric classification, published in Enziklopaedie der Mathematische Wissenschaft\footnote{\tiny See \cite{54} \it Kontinuierliche geometrische Gruppen; die Gruppentheorie als geometrisches Einteilungsprinzip, }, is a natural example of this fact. It is noteworthy that this contribution of Fano was taken up by Elie Cartan, which around 1914  translated it, writing also some additional parts, for the French edition of Enzyclopaedie\footnote {\tiny E. Cartan \it Oeuvres compl\'etes \rm volumes I-VI, reprinted by CNRS (1984)}. \par Fano's contributions  to Lie theory, in particular on finite representations of the complex special linear group and the finite generation of its ring of invariant polynomials,  are described by Armand Borel in his historical essay on Lie groups \footnote{\tiny A. Borel \it Essays in the history of Lie groups and algebraic groups\rm, History of Mathematics v. 21, AMS Providence RI (2001)}. Fano's contributions also appear in the section dedicated to him in the book of Thomas Hawkins on history of Lie groups\footnote{\tiny  T. Hawkins \it The emergence of Lie groups. An essay in history of mathematics (1869-1926)\rm Springer  (2000)} \par
It is therefore obvious that the new interplay between geometry and physics, appearing in the late 19th and early 20th centuries, was perfectly clear to him, as well as it is clear his interest on it. This can be appreciated in many ways, for instance from his book on non euclidean geometries, with the explicit subtitle 'Geometric introduction to general Relativity' \footnote{\tiny \cite{113}: \it Geometria non euclidea, (introduzione geometrica alla teoria della relativit\`a generale).}  or from conferences and articles due to him. \par However physics had to become even more important
 in Fano family. Two are the sons of Gino Fano:  one is Robert Fano, born in 1917 and presently Professor Emeritus of Electrical Engineering and Computer Science at Massachussetts Institute of Technology. \par The other one
 was Ugo Fano (1912-2001). Ugo Fano, Enrico Fermi award 1996 and Fermi's student and collaborator, is well known for his fundamental contributions to atomic physics.  Let us quote some memories about him and his father from New York Times: \par { \it  'Ugo Fano was born in Turin, Italy, on July 28, 1912. His father, Gino Fano, was a mathematician at the University of Turin, and he kept the young Ugo apprised of the great discoveries in physics and mathematics that were taking place around the world. In a memoir, Ugo Fano said he recalled being introduced to Niels Bohr's new atomic theory at the dinner table when he was 12' \footnote{\tiny J. Glanz \it Ugo Fano is dead at
 88; Physicist linked to Fermi \rm New York Times, February 15 2001.}. }\par
The latter remarks put in evidence the broadness of the scientific views of Fano. This is quite typical of great men of Science born at his times. However it is perfectly combined with his outstanding mathematical talent and the deep relations he had with the scientific milieu of his times. \par To this respect,  before of passing to the main stream of his algebro-geometric interests, the picture could be expanded in several directions, in order  to appropriately describe the presence and the strong influence of Fano in the mathematical and scientific debate of his times.  \par Let us introduce just one more example of this fact. Here we refer to the intense debate on foundations of Geometry in the years preceding Hilbert's Grundlangen der Geometrie. \par
As is well known this debate was of special importance in Italy. The University of Torino became one of the main centers of it, due to the presence of Corrado Segre and Giuseppe Peano: two leading scholars with their own points of view on foundations of Geometry. \par Many geometers contributed to this subject and the young Fano was among them. His ideas on projective geometry in arbitrary dimension are near to Segre's ideas and to $n$-dimensional geometry of Bertini and Veronese. \par They are neatly expressed in his paper of 1892 on foundations of projective geometry, where a famous model of projective plane, named today the Fano plane, is in particular constructed. \par It is not our task to enter in the historical studies
on those years, it is however appropriate for our description mentioning the relevance of this contribution given by Fano. \par Indeed, about the Italian debate on foundations of Geometry, it is not rare to read comments in the same spirit 
of the following one, due to Jeremy Gray \footnote{\tiny  J. Gray \it The Foundations of Projective Geometry in Italy \rm Undergraduate Mathematics Series p. 269-279 Springer (2010)}. He is essentially reporting Hans Freduenthal's point of view: \par
{\it  'When the distinguished mathematician and historian of mathematics Hans Freudenthal analysed Hilbert's Grundlagen he argued that the link between reality and geometry appears to be severed for the first time in Hilbert's work. However, he discovered that Hilbert had been preceeded by the Italian mathematician Gino Fano in 1892. Recent historians of mathematics have shown that, in Italy at least, Fano's point of view on the nature of geometrical entities had been a generally accepted theory for at least a decade, but it was not in fact axiomatic in Hilbert's manner; other Italian mathematicians were, however, ahead of Hilbert in this regard. How did this come about, what did they do, and why they lose out?'
} \par
 We have stressed enough the wideness of Fano's views and his strength in the general field of geometry.  Nevertheless Fano was mainly, though definitely not only, an algebraic geometer. \par In particular he was immersed in the main stream of classical Italian school of algebraic geometry. This can be seen  in many ways: for instance from the central role played by birational geometry in his views or  from the attention he payed to classification  of algebraic varieties. \par Another example of this fact is certainly represented by his genuine and strong intuition in understanding the natural geometric life, so to say, of a projective variety living in an $n$-dimensional projective space.  However, in the taste of Fano as an algebraic geometer and in his choices of investigation, something more  is present which goes further. \par Taking our sentence with a grain of salt, we can say that, often, he opened the way to problems, or to important geometric objects, which were quite  unexplored or still untouched: \par birational geometry in higher dimensions, algebraic groups of Cremona transformations, cubic threefolds and their irrationality, K3 and Enriques surfaces with special automorphisms, Fano threefolds, of course, and more.  
 \newpage
 \section{On Fano's research and scientific heritage }
In what follows we enumerate the main research topics Fano addressed. Each of them is accompanied by a short note, whose purpose is to outline Fano's contribution and the periods of his work on the subject.  \par Each note consists of free and non systematic comments: we usually single out just one representative paper and use it to very briefly describe the subject and Fano's results.  At the end of any subject we quote some relevant papers of Fano on it. \par   The subjects of the list follow, approximately, the chronological order of the scientific events in the life of Fano. If possible our aim is to capture, by a unique rapid sketch, much of Fano's scientific image.   
 \begin{itemize}
\item[$\circ$] \underline {\it  Foundations of Geometry} 
\end{itemize} { \small   
This subject was specially cultivated by the very young Fano.  As already remarked, his contribution is situated in a fertile area of research of those years, in particular in Italy\footnote{\tiny To the  mentioned names of G. Peano, C. Segre and G. Veronese, we should add at least those of F.Amodeo, F. Enriques, M. Pieri, G. Vailati and many others.}. His paper \cite{2}\footnote{\tiny  \it Sui postulati fondamentali della geometria proiettiva in uno spazio lineare a un numero qualunque  di dimensioni}  is of relevant importance. \par At that time $n$-dimensional projective spaces were currently in use, and investigated, via projective coordinates. However there weren't established foundations of this notion by a system of independent postulates,  allowing in particular the use of coordinates. This important problem was raised by Corrado Segre\footnote{\tiny See \it I Quaderni di Corrado Segre, \rm CD-ROM edited by L. Giacardi, Dipartimento di matematica, Universit\`a di Torino, 2002} to his  students, among them Fano, during his course 'Geometria sopra un ente algebrico' in the academic year 1890-91.  \par In \cite{2} Fano constructs an adequate system of independent postulates for a projective $n$-space\footnote{\tiny Notably he says that points in this space are `\it entities which, for brevity, we shall call points, independently, however, from their nature'\rm, \cite{2} p. 109.} and, to confirm their independence, he comes to discover examples of finite projective spaces. By an ingenious construction he constructs in this way the projective plane over the field $\mathbb Z/2 \mathbb Z$, still  named the Fano plane. Notably these examples represent early sources for the field of finite geometry. Fano's ideas on foundations of geometry  are expressed in several papers, articles and letters.}  \par
{ \em Some contributions along the years\rm: \cite{1, 2, 10, 11, 36, 55, 59, 68, 82, 106, 107, 111, 118, 119, 123}  \medskip
 \begin{itemize}
\item[$\circ$] \underline{ \it  Algebraic curves} 
\end{itemize}
{\small The theory of algebraic curves, whose place is a central one in the interests of the classical Italian school, shows a minor weight in the interests of Fano. Curves are everywhere in his
papers, but rarely these are concerned with this theory. \par A remarkable exception, influenced by the scientific exchange with Castelnuovo and Segre, is the paper \cite{4}. Actually this is part of 
'Dissertazione di Laurea'\footnote{\tiny The thesis to be submitted by a university student for obtaining the degree.}, written by Fano under the direction of Segre and published in 1894. \par
It deals with curves $C$ of degree $d$ in a projective space $\mathbb P^r$. Two years before Castelnuovo had proved that the genus of $C$ is  bounded and that the maximal genus $p(d)$ is reached by curves in a rational normal scroll or in the Veronese surface\footnote{\tiny G. Castelnuovo  \it Ricerche di Geometria sulle curve algebriche \rm 24, Atti Acc. Scienze Torino (1890)}. \par Fano applies
the same methods to characterize curves of genus near to $p(d)$. In particular he studies the linear system of quadrics through a curve $C$ of genus $p(d)-k$. Together with further results he deduces that, for $k = 1$,  $C$ is contained in a surface whose hyperplane sections are rational or elliptic. From this property the classification of curves $C$ in genus $p(d) - 1$ follows.}
\par
 { \em Some contributions along the years\rm: \cite{2, 98, 126} \medskip
 \begin{itemize}
\item[$\circ$] \underline { \it  Continuous groups in projective and birational geometry}
\end{itemize}  
{\small Due to the activity of Klein and Lie, transformation groups were a central theme on the mathematical scene during the early years of Fano.  He was well exposed to these ideas due to
the scientific ties with Klein we mentioned. As a consequence, the theme of continuous groups of birational transformations represents a fundamental part of Fano's work.  He was specially active on it  in the five years before 1900. Although other contributions are possibly better known, this one consists of notes, memoirs and essays, for more than 300 pages, of maximal interest. \par In the spirit of Erlangen program a leading project on the field was to understand projective varieties, in the complex projective space $\mathbb P^n$, which are invariant by the action of any given Lie subgroup $G$ of $PGL(n+1, \mathbf C)$. This research was started by Klein and Lie in $\mathbb P^2$ and $\mathbb P^3$. They were joined by Enriques and others in the same period. The classification came to completion  before 1900\footnote{\tiny See: (1) F. Enriques `\it Le superficie con infinite trasformazioni proiettive in s\`e stesse' \rm Atti R. Ist. Veneto 51 (1893) 1590-1635.  (2) S. Lie `\it Bestimmung aller Flaechen, die eine continuierliche Schaar von projektiven Transformationen gestatten \rm in 'Gesammelte Abhandlungen' Teubner, Leipzig (1922). \ Various contributions are due to Fano. }.   \par
In two memoirs, \cite{17, 20}, Fano extends this classification to $\mathbb P^4$. It is not a kind of trivial step, he uses new ideas from geometry and representation theory.  Actually, when $G$ is a non integrable\footnote{\tiny In the classical language of Lie's theory \it integrable group \rm precisely means solvable group, see Thomas Hawkins loc. cit. p. 91. Since this word
often appears in the papers of Fano and in their titles, we have chosen to keep it.} continuous group, his method turns out to be the search in $\frak {sl}(5, \mathbf C)$ of Lie subalgebras isomorphic to $\frak {sl}(2, \mathbf C)$. He applies Study's complete reducibility theorem for $\frak{sl}(2, \mathbf C)$ and also gives a geometric proof of it. Then he can construct the beautiful series of $G$-invariant  hypersurfaces. \par
A third fundamental memoir is a joint paper of Enriques and Fano, \cite{25}. It is of primary importance in the study of the Cremona group $Bir (\mathbb P^3)$, the group of Cremona transformations, that is birational automorphisms, of $\mathbb P^3$. In it all the complex Lie subgroups of  $Bir(\mathbb P^3)$ are classified. \par The authors show that, up to two exceptions, each of them is conjugated to a known subgroup of one of the following continuous groups: $PGL(4, \mathbf C)$, the group of conformal transformations, the two types of groups defined by de Jonqui\'eres space transformations. The exceptions are two 3-dimensional groups respectively related to the octahedron and the hicosahedron. \par
Finally let us recall that de Jonqui\'eres transformations of $\mathbb P^3$ generalize those of $\mathbb P^2$ with the same name. They are of two types according to what they fix: a pencil of planes or a star of lines. In a fourth memoir Fano  classifies all continuous groups of them, \cite{27}. This complements the joint memoir with Enriques. \par  Continuous groups of Cremona transformations were object of a communication by Fano at the first International Congress of Mathematicians, at Zurich in 1897. } \par
 {\em Some contributions along the years\rm: \cite{12, 17, 18, 20, 21, 22, 24, 25, 26, 27, 28, 30, 31, 32, 38, 54}} \medskip
 \begin{itemize} \item[$\circ$] \underline{ \it  Algebraic varieties defined by linear differential equations}
\end{itemize}
{\small This topic represents a very important contribution due to  Fano. It is coherent with the previous one on continuous groups and it was developed 
in the same years. In particular, it can offer to the reader a more precise and comprehensive focus on the real features and complexity of Fano's work  in its historical frame. \par
This contribution reflects once more the influence of the circle of  ideas arising from the impulse of Lie's work, that is, the foundations of Galois theory of linear differential equations\footnote{\tiny Cfr. M. van der Put, M. F. Singer \it Galois theory of linear differential equations \rm  Grundlehren Math. Wiss. vol. 328 Springer Berlin (2003)}. From this point of view Fano was developing the classical theory of Picard-Vessiot, as he explains in \cite{34} and in other papers. \par Fano applies this theory to the remarkable case where there exists a non zero homogeneous polynomial
$F \in \mathbf C[T_1 \dots T_n]$ which is identically zero on a set of $n$ independent solutions $y_1 \dots y_n$ of a given linear differential equation $L = 0$.  \par The Galois group of $L = 0$ acts
linearly on the projective subvariety $M \subset \mathbb P^{n-1}$, defined by all the polynomials $F$. Moreover, the maximal continuous group of projective transformations fixing $M$ yields most of the informations about the integration of $L = 0$.  This relates the study of $L$ to the classification of projective varieties which are invariant by a continuous subgroup of $\Aut  \mathbb P^{n-1}$.\par 
Most of Fano's results on this topic converge to a fundamental paper, published by him in Mathematische Annalen, \cite{39}. Notably this paper, as well as his entire work in the field, were reconsidered by Michael Singer. Singer offers a very clear description of Fano's contribution in his paper  \it Algebraic Relations Among Solutions of Linear Differential Equations: Fano's Theorem \rm \footnote{\tiny M. Singer in American Journal of Math. 110 (1988) 115-143.}, published in 1988. Let us freely summarize from it:  \par Assume $L = 0$ is homogeneous with complex rational functions as coefficients. Then, with the previous notations, consider the following statement:  \par (*) \it If there exists a non zero $F \in \mathbf C[T_1 \dots T_n]$ such that $F(y_1 \dots y_n) = 0$, then all solutions of $L = 0$ can be expressed as algebraic combinations of solutions of linear differential equations of order $\leq n-1$. \rm \par
A major Fano's contribution is the theorem that (*) is true for $n \leq 5$. He also gave partial positive results for $n \geq 6$. Due to these and later results\footnote{\tiny Cfr. D.V. Chudnovsky, G.V. Chudnovsky \it The wronskian formalism for linear differential equations and Pade approximation \rm Advances in Math. 53 (1984) 28-54}, (*) was expected
to be true for $n \geq 6$. Finally Singer showed that (*) is true for $n = 6$ and false for $n \geq 7$. \par Fano's approach to (*)  is very explicit and beautiful. He uses group theoretic and projective geometric techniques to detect the possible algebraic relations satisfied by $y_1, \dots, y_n$ for $n \leq 5$. Then he relies on his deep knowledge of projective varieties in $\mathbb P^{n-1}$, $n \leq 4$, which are invariant by a continuous subgroup of $\Aut  \mathbb P^{n-1}$. } \par
 { \em Some contributions along the years\rm: \cite{7, 13, 14, 15, 16, 33, 34, 35, 37, 39} } \medskip
 \begin{itemize}
\item[$\circ$] \underline {\it Line geometry}  
\end{itemize}
{\small Line geometry is the study of the Grassmannian $G(1,\mathbb P^n)$ of lines of a projective space $\mathbb P^n$. Around 1900 Fano dedicated his best to this subject, aiming to
continue the classification of congruences of lines started by Kummer over the complex field\footnote{\tiny  E. Kummer    `\it Ueber algebraische Strahlen systeme, insbesondere ueber di der ersten und zweiten Ordnung \rm' Berliner Abh. 1866}. A congruence of lines is just a surface $S$ in $G(1, \mathbb P^3)$. Its cohomology class defines a pair of
integers $(a,b)$, known as the type of $S$. \par At that time, after a long work, the classification of congruences of lines was known for $a \leq 2$. Its extension to the case $a = 3$ was essentially realized
by Fano in the very rich paper 'Congruenze  di rette del 3$^o$ ordine prive di linea singolare', \cite{39}, and in some related papers. \par
Notably Fano abandons, at least partially, the study in $\mathbb P^3$ of families of lines. His modern geometric approach is openly declared \footnote{\tiny ` \it Studio di alcuni sistemi di rette considerati come superficie dello spazio a cinque dimensioni' \rm is the title of \cite{3}.}: to consider and study these families as subvarieties of the appropriate space 
for line geometry. This means the Pluecker embedding of $G(1, \mathbb P^3)$ that is the Klein quadric in $\mathbb P^5$ . \par Our distance from Fano's results
is not big: problems on  congruences of lines and their classification are still widely open.  \par
Line geometry of $\mathbb P^n$ is moreover considered by Fano in two remarkable cases. They are related to cubic hypersurfaces in $\mathbb P^4$, a crucial topic of his investigations.
\par For $n = 4$ he studies the family $F(V)$ of all lines contained in a smooth cubic hypersurface $V \subset \mathbb P^4$, \cite{45}. 
$F(V)$ is a surface, known  as the Fano surface of $V$.\par  For $n = 5$ he studies the general  3-dimensional linear section of the Pluecker embedding of $G(1, \mathbb P^5)$, proving
 that it is birational  to a general cubic $V$, \cite{100}. } \par
 {\em Some contributions along the years\rm:  \cite{3, 5, 19, 39, 41, 43, 45, 49, 78, 99, 100} \medskip
 \begin{itemize}
\item[$\circ$] \underline {\it Non rationality and birational geometry in dimension 3}  
\end{itemize}
{\small In 1876 Lueroth proved that a unirational complex curve is also rational. The extension of this theorem to higher dimensions became few years later the Lueroth problem.
In 1892 this was positively answered in dimension two by a brilliant theorem of Castelnuovo\footnote{\tiny See: G. Castelnuovo `\it Sulla razionalit\'a delle involuzioni piane \rm' 1894, Math. Annalen 44 (1894) 125-155}. As is well known the negative answers to Lueroth problem in dimension $\geq 3$ are a central
episode of XX century algebraic geometry. \par In the paper \cite{56} Fano gives a proof of non rationality for a general quartic hypersurface in $\mathbb P^4$ and for a general quadro-cubic
complete intersection in $\mathbb P^5$. Though his proof is not satisfactory from
a modern point of view, the main ideas behind it are the correct ones. \par They were fruitfully cultivated by the Russian school of Manin and Iskovskih. The non rationality
of the quartic then followed by Iskovskih and Manin in 1971\footnote{\tiny For a general picture see the volume '`\it Algebraic Geometry V - Fano Varieties \rm'  by V. Iskovskikh and Yu. Prokhorov in Encycl. of Math. Sciences. v. 47 (1991) Springer Berlin}. Later the same method of birational rigidity was successfully applied by Pukhlikov to the other threefold\footnote{\tiny A. Pukhlikov  `\it Maximal singularities on the Fano variety $V^3_6$\rm'  Moscow Univ. Math. Bull. 44 (1989)  In this case, however, Beauville had already proven the non rationality via the intermediate jacobian.   70-75. \ A. Beauville `\it Vari\`et\`es de Prym et jacobiennes interm\`ediaires \rm' Ann. Sci. Ec. Norm. «Sup. 10 (1977) 309-391, thm. 5.5.}.  \par
Fano's approach to this problem is today encoded in the notion of birational rigidity and relies on a technique for the factorization of a birational map which is known as Noether-Fano inequality
\footnote{\tiny See the mentioned book 'Algebraic Geometry V - Fano Varieties', Lemma 9.9.1 }. The prototype of it is the inequality,
due to Noether, used for proving that every Cremona transformation of the plane factors through quadratic and linear transformations. \par
In the related paper \cite{65} Fano improves his program in order to prove the non rationality of  a threefold $V$ as above. He proposes a new, very modern, idea:  to study the linear systems of K3
surfaces existing on $V$ to deduce its birational rigidity and hence its irrationality. \par On the other hand, since $V$ is in the series of Fano threefolds, the question of its unirationality of is natural.
If $V$ is a general quadro-cubic complete intersection this was proved by Enriques in 1912. Hence $V$ represents in this case a negative answer to Lueroth problem. Still today the unirationality of a general quartic of $\mathbb P^4$ is instead an outstanding open problem. \par
It is due to mention in this section the interest of Fano for developing further techniques of 3-dimensional birational geometry . This is a motivation for his work on the classical theory of contact birational transformations of the plane, which is of projective-differential geometric flavor. \par Let $\mathbb P \Omega^1_{\mathbb P^2}$ be the projectivized cotangent bundle of $\mathbb P^2$. These transformations can be defined as birational automorphisms $a$ of $\mathbb P \Omega^1_{\mathbb P^2}$  such that $p \cdot a = b \cdot p$, where $p: \mathbb P\Omega^1_{\mathbb P^2} \to \mathbb P^2$ is the natural projection and $b \in Bir(\mathbb P^2$).  Since $\mathbb P \Omega^1_{\mathbb P^2}$ is rational, $a$ is birationally equivalent to an element of $Bir (\mathbb P^3)$.  
 \par To this, quite unexplored, topic is dedicated one of the two communications of Fano to the International Congress of Mathematicians of 1928.}  \par
{ \em Some contributions along the years \rm: \cite{27, 30, 38, 56, 65, 66, 70, 94, 95, 102, 108, 129, 130, 137, 143}} \medskip
\begin{itemize}
\item[$\circ$] \underline{ \it  Cubic threefolds} 
\end{itemize}
{\small Cubic hypersurfaces in the complex 4-dimensional projective space are usually called cubic threefolds. They are definitely associated to Fano's biography, since he was studying cubic threefolds  systemically for a very long period, with a special attention towards the problem of their irrationality. \par
Since the unirationality of a smooth cubic threefold $V$ was known\footnote{\tiny Virgil Snyder in `\it The problem of the cubic variety in $S_4$\rm', (Bull. AMS 35 (1929) 607-642) says that this property was known to Max Noether}, proving its irrationality was sufficient for Fano to negatively answer Lueroth problem. This was of course one of his motivations for studying cubic threefolds.
 In the late survey \cite{143} of 1950 Fano describes his approach to the irrationality of $V$ and many other aspects of his work on cubic threefolds. \par This begun around 1904: in this year he  publishes four papers on cubic threefolds where various 
 fundamental properties are established. He shows that the Picard group of a smooth $V$ is generated by the hyperplane class, so that every
 surface in $V$ is complete intersection. Also, he studies in detail the dual hypersurface and more properties of a general of $V$. \par
 Another cornerstone is his study of the family of lines contained in $V$. This is a surface $F(V) \subset G(1, \mathbb P^4)$, known as the Fano surface of $V$. As is well known it plays a crucial role in the proof of the irrationality of $V$ due to Clemens and Griffiths via the method of the intermediate jacobian\footnote{\tiny H. Clemens, P. Griffiths ` \it The intermediate jacobian of the cubic threefold' \rm  Annals of Math. 95 (1972) 281-356}. It is worth mentioning a geometric investigation by Fano which is quite connected to this method. It is the study of plane quartics which are everywhere tangent to a given plane quintic, \cite{126}. \par Actually any line in $V$ defines one of these families of quartics. Moreover the theta divisor  $\Theta$ of the intermediate jacobian of $V$ can be reconstructed from the corresponding family of even theta characteristics of these quartics. In modern terms he considers families of even spin curves of genus $3$ which are birational to the quotient of $\Theta$ by $-1$ multiplication. \par
 Soon Fano realized that the methods of birational geometry he could use at his time were not satisfactory in the case of $V$. This opened the way to him for studying more properties
 of $V$ and a wider series of threefold relatively similar to $V$. \par Today these are named Fano threefolds. The 3-dimensional linear section $V'$ of the Pluecker embedding of the Grassmannian $G(1, \mathbb P^5)$ is one of them.
 Fano proves the beautiful result that a general $V$ is birational to some $V'$ and conversely.  Then he gives a proof of the irrationality of $V'$ via the methods mentioned in the previous subject. }}\par
   {\em Some contributions during the years\rm:  \cite{45, 46, 47, 48, 65, 70, 99, 100, 126, 131, 133, 136, 143}} \medskip \begin{itemize}
\item[$\circ$] \underline{ \it  Fano threefolds}
\end{itemize}
{ \small A smooth projective variety $V$ is a Fano variety if its first Chern class is ample. The definition suitably extends to a singular $V$.  In contemporary algebraic geometry Fano
varieties are among the fundamental elements for the birational classification of algebraic varieties of any dimension. They bear the name of Fano because of his pioneering work in the classification of Fano threefolds i.e.
3-dimensional Fano varieties. \par
Starting from the study of the  irrationality problem for some of these threefolds, like for instance cubic threefolds, the activity of Fano on this subject developed along forty years.  After his death the investigations on Fano threefolds became a  
fundamental direction of algebraic geometry, culminating in the complete classification of them, based on modern methods and on Mori theory \footnote{\tiny See again Algebraic Geometry V - Fano Varieties}. Nonetheless Fano's previous contributions to this 
classification are everywhere present and recognized as invaluable. 
In 1928 his first results were communicated to the International Congress of Mathematicians held in Bologna.\par  In two memoirs of 1937 he gave a bound on the degree of Fano threefolds, thus proving  that they are distributed in finitely many families. He also proved most of the possible rationality results for Fano threefolds. Later he recollected most of his work in the memoir 'Nuove ricerche sulle variet\'a algebriche a tre dimensioni a curve sezioni canoniche', published in 1948. \par In it one finds a deep geometric description of all the families he had previously considered. Another part is dedicated to  Lueroth problem and to the irrationality problem for the members of some of these families. \par
Fabio Conforto, an emerging algebraic geometer of that time, made a lucid and fascinated analysis of this memoir in a long review\footnote{\tiny F. Conforto Mathematical Reviews MR0038100.}. It is interesting to read in it the profound value he attributes to this work, even if he remarks some known limits reached by the classical methods in use. } \par
   {\em Some contributions along the years \rm: \cite{46, 56, 99, 100, 103, 108, 114, 115, 116, 120, 121, 129, 130, 135, 136, 139, 140, 143}} \medskip
 \begin{itemize}
\item[$\circ$] \underline {\it   Fano-Enriques threefolds}
\end{itemize}
{\small Let $V \subset \mathbb P^n$ be a Fano threefold. One says that $V$ is anticanonically embedded if its first Chern class is the class of its hyperplane sections. In this case a general hyperplane section of $V$ is
a K3 surface\footnote{\tiny A K3 surface is a compact complex surface which is simply connected and such that its first Chern class is zero.}. The interest of Fano for K3 surfaces has certainly one motivation in his thorough
investigations on linear systems of K3 surfaces on a Fano threefold. Another motivation, we rediscover in the next section, is his ancient curiosity for automorphisms groups in algebraic geometry. \par
Enriques surfaces are very related to K3 surfaces: an Enriques surface is the quotient of a K3 surface endowed with a fixed point free involution. Now, in the world of Fano threefolds, it appears that some of them admit an involution with
finitely many fixed points, namely eight. \par Taking the quotient of such a threefold by such an involution, one obtains a threefold $V$ which is a Fano threefold again. It is of special type. In particular, a birational projective model of $V$ is a threefold $V'$ whose hyperplane sections  are Enriques surfaces. Today we say that $V'$ is an Enriques-Fano threefold, just because its hyperplane sections are Enriques surfaces. \par
Enriques-Fano threefolds are a kind of geometric jewels in the landscape we are describing.  Fano met the study of them in the memoir, \cite{122}. His results parallels those of Godeaux\footnote{\tiny L. Godeaux ` \it Sur les variet\'es alg\'ebriques \`a trois dimensions dont les sections sont des surfaces de genre z\'ero e de bigenre \rm' un' Bull. Acad. Belgique, Cl. des Sci. 14 (1933) 134-140.} in the same years. \par  He proves  that $V'$ is birational to a Fano threefold. Then he gives a beautiful geometric classification 
of these threefolds $V'$. This is however based on a restrictive hypothesis on the eight singular points  present in $V'$.  
Fano-Enriques threefolds were essentially rediscovered in the eighties\footnote{\tiny  A. Conte J.P. Murre ` \it Algebraic varieties of dimension three whose hyperplane sections are Enriques surfaces\rm ' Ann. Sc. Norm. Sup. (4) 12  (1985) 43-80}.  As in the case of Fano threefolds, the work of Fano, though not bringing a complete classification, strongly influenced modern investigations, which became flourishing on this topic.\footnote{\tiny See for instance: Y. Prokhorov \it On Fano-Enriques varieties \rm Sb. Math. 198 (2007) 559-574 and A. Knutsen, A. Lopez,  R. Munoz \it On the extendability of projectives surfaces and a genus bound for Enriques-Fano threefolds \rm J. Diff. Geom. 88 (2011) 485-518}}  \par
{\em Some contributions along the years: \cite{122,134}} \medskip
\begin{itemize}
\item[$\circ$] \underline {\it   K3 or Enriques surfaces and their automorphisms}
\end{itemize} {\small The automorphisms groups of K3 or Enriques surfaces are discrete, trivial for a general K3. The work of Fano on these groups of automorphisms is linked to quartic surfaces and also to his investigations on line geometry.
   \par In 1906 Enriques proved the following theorem: every surface having a non finite, discrete group of automorphisms is either elliptic or a K3 surface. In the same year Fano constructs the first example of a non elliptic K3 surface of this type. \par It is a quartic surface of Picard number two containing a smooth sextic of genus two, \cite{51}]. 
His conjectural remark on the existence of infinite series of families of examples, distinguished by numerical characters, is completely confirmed. Later he came to enlarge the picture constructing several families, see \cite{71, 72, 73,74, 75, 132, 141}. \par
But a beautiful jewel was constructed by him when studying automorphisms of Enriques surfaces. Enriques surfaces are elliptic. Moreover, as observed by Enriques,  the group of automorphisms is not finite for a general  Enriques surface. \par In particular it  could appear almost unbelievable, or at least not obvious, that this group could become finite for special elements in the family of all Enriques surfaces{\footnote{\tiny Cfr. V. Nikulin 'On algebraic varieties with finite polyhedral Mori cone' in The Fano Conference, Editors A. Collino, A. Conte, M. Marchisio.  printed by University\'a of Torino, Dipartimento di Matematica Torino (2004) 573-589, section 3.}. In 1910 Fano produced an example of this type, \cite{60}.  \par The example is a very special case of Reye congruence. A Reye congruence $S$ is a special  Enriques surface embedded in the Grassmannian $G(1, \mathbb P^3)$. Its Pluecker embedding is known as an example of Fano model of an Enriques surface\footnote{\tiny  F. Cossec, I. Dolgachev-Cossec `\it Enriques Surfaces I\rm' Birkha\"user, Boston (1989) p. 279 and I. Dolgachev `\it A brief introduction to Enriques surfaces\rm', conference in honor of S. Mukai http://www.math.lsa.umich.edu/~idolga/Kyoto13.pdf}. \par Roughly speaking the very special $S$ considered has an elliptic pencil which is invariant by $Aut \ S$ and moreover contains a very special element $F$. $F$ is a suitable configuration of rational curves, which is invariant by $Aut \ S$. The existence of such an $F$ forces $Aut \ S$ to be finite. \par Enriques surfaces with a finite group of automorphisms disappeared from the scene after the isolated construction of Fano. One had to wait the eighties of XX century for their rediscovery. This was done by Igor Dolgachev in 1983. Further results and the classification followed\footnote{\tiny See: I. Dolgachev `\it On automorphisms of Enriques surfaces\rm' Invent. Math. 76 (1984) 163-177,  W. Barth, C. Peters `\it Automorphisms of Enriques surfaces\rm',  73 Invent. Math. (1983) 383-411}. Later Fano's example was rediscovered by Dolgachev himself and discussed from a modern point of view.}}\par {\em Some contributions along the years: \cite{26, 51, 60, 71, 72, 73, 74, 75, 122, 132, 134, 141}} }
 \newpage
  \section{Epilogue}
\rm  For most of his life Gino Fano was professor  of projective and descriptive geometry at the University of Torino. More precisely he was professor there since 1901 until the year 1938. \par
 In this year he lost his job due to the racial laws against jews established by the fascist regime. Then he had to leave the country. He was reintegrated in his job in 1946, after the end of World War II. \par 
 In contrast to these obscure years, it is perhaps the right moment for mentioning the deep and positive influence of Fano, as a person and as a teacher, on so many persons and brilliant students, as well as mentioning the sign he left 
 during the long years spent as a professor. \par In the year 1950 a group of persons, made by former students 
 and friends, organized in Torino a series of lectures by him as a homecoming celebration in his honor.  Reading the list of these persons, in the volume of this event, is impressive, and instructive of the above mentioned influence of Fano. \par 
 In the list, together with other famous mathematicians, one can see the name of Beniamino Segre. His words
 can serve as an epilogue for a long, great life in geometry. Perhaps they can serve for an entire geometric age, in the years of transition from older to newer algebraic geometry: \bigskip \bigskip \par \par  '... \it con riconoscenza ed ammirazione al Professor Gino Fano, le cui smaglianti lezioni di Geometria proiettiva - nel lontano 1919-20 - ebbero su di me un effetto decisivo, attraendomi irresistibilmente verso gli studi geometrici'. \footnote {\tiny B. Segre `\it Intorno agli spazi lineari situati sulle quadriche di un iperspazio\rm' Rend. Sem. Mat. Univ. e Polit. di Torino 9, 1949-50, 137-144.} \rm
 \newpage
  \section{Gino Fano: 1871 - 1952} 
  \small  \medskip \par \noindent
$\bf 1871$ \it Gino Fano\footnote{\tiny For further bio-bibliographical informations we quote and suggest the very interesting monography  \it Sull' apporto evolutivo dei matematici ebrei mantovani nella nascente nazione italiana,\rm   by A. Janovitz and F. Mercanti, Monografie di EIRIS(2008), rivista on line www.eiris.it} was born on January 5, 1871, in Mantova, Italy.\footnote{ \tiny Gino Fano is the son of Ugo Fano and Angelica Fano. His family was a wealthy Jewish family rooted in Mantova. His father volunteered with Garibaldi. His mother was quite engaged in defense of Italian Culture and member of 'Dante Alighieri' Society.} \medskip \par \noindent
$\bf 1871$ In the same year and day Federigo Enriques was born in Livorno.  \medskip \par \noindent
$\bf 1871$ Rome is the new capital city of the Kingdom of Italy. \medskip \par \noindent
$\bf 1880-87$ Fano has a brilliant career as a student in various schools.\footnote{\tiny 1880-83: Fano studied at Liceo-ginnasio Virgilio in Mantova. Obeying his father, he started a military career.  Then he left it for the section of Physics and Mathematics of R. Istituto A. Pitentino in Mantova, where he brilliantly obtained his diploma and a special prize as an excellent student.} \medskip \par \noindent
$\bf 1888-92$ Student at University of Torino, Fano is pupil and friend of Segre and Castelnuovo. Suggested by Segre, he translates Klein's Erlangen Program. \medskip \par \noindent
$\bf 1892$ June 22: Fano is proclaimed Dottore in Matematica, his advisor is Segre.\footnote{\tiny The committee: G. Basso, E. d'Ovidio, N. Jadanza, G. Peano, M. Pieri and C. Segre. The grade given to Fano was maximal: 90/90 cum laude.} The contents of his Dissertazione di laurea become a published paper, \cite{4}.
\medskip \par \noindent 
$\bf 1892-93$ Assistant of E. d' Ovidio at the University of Torino. First contributions to foundations of geometry and to line geometry.
\medskip \par \noindent
$\bf 1893-94$ Academic year spent in Goettingen under the invitation of Felix Klein.  
\medskip \par \noindent
$\bf 1894-95$ Assistant of G. Castelnuovo in Rome. Contacts with Luigi Cremona and various mathematicians in Rome.
\medskip \par \noindent
$\bf 1892-97$ Several contributions to:  continuous groups and their invariant projective varieties, geometric properties of
linear differential equations, line geometry. 
\medskip \par \noindent
$\bf1897$ First International Congress of Mathematicians, Zurich: communication by Fano on continuous groups of Cremona transformations.
\medskip \par \noindent
$\bf 1899$ Fano wins a position of professor on the chair of algebra and analytic geometry of the University of Messina.
 \medskip \par \noindent
$\bf 1899$ A position of professor in Goettingen, starting from December 1st 1899, is offered by Klein to Fano. However he will decline this offer.\footnote{\tiny   According to Fano's son Ugo, his father '...did not want to be Germanized'.}\medskip \par \noindent
$\bf 1899 $ Hilbert's Grundlagen der Geometrie. \medskip \par \noindent
$\bf 1899 $ Social conflicts and industrialization in Italy.  July 11: opening of FIAT, Fabbrica Italiana di Automobili - Torino. 8 cars are produced in 1899. \medskip \par \noindent $\bf 1900$ Second International Congress of Mathematicians: Hilbert's problems. \medskip \par \noindent
$\bf 1901$ Fano becomes professor in the University of Torino, on the chair of projective and descriptive geometry. \medskip \par \noindent
$\bf 1897-03$ Memoirs on continuous groups of Cremona transformations. Memoirs on order $3$ congruences of lines. Notes and long papers on linear differential equations whose
solutions satisfy algebraic relations.\medskip \par \noindent
$\bf 1904$ First papers on the cubic threefold: a series of four papers. One of them describes the Fano surface of lines of a cubic threefold. \medskip \par \noindent
$\bf 1904$ 11-th Summer Meeting of American Mathematical Society, St. Louis: Henry Poincar\'e and Gino Fano are special invited guests of the Meeting.
\medskip \par \noindent
$\bf 1905$ Albert Einstein's Annus Mirabilis\footnote{\tiny Four revolutionary papers published by Einstein in 1905 on Photoelectric effect, Brownian motion, Special relativity, Mass-energy equivalence.}\medskip \par \noindent
 $\bf 1905$ Fano begins his cooperation, as a teacher and as an organizer, to the Evening School for Women Workers of Torino.\footnote{\tiny For his action in favor of the School and Public Education, Fano will get the golden medal of 'Benemerito della Pubblica Istruzione' in 1928} See his report on the period 1909-10.
 \medskip \par \noindent 
$\bf 1907$ Two important contributions of Fano to Enzyclopaedie der Mathematische Wissenschaft, see \cite{53, 54}
\medskip \par \noindent
$\bf 1907-11$ First papers on K3 or Enriques surfaces, in particular on their automorphisms. First papers on the rationality problem for some Fano threefolds.
\medskip \par \noindent
  $\bf 1911$ Torino: Gino Fano marries Rosa Cassin (1880-1956).They will have two sons, Ugo (1912-2001) and Robert (1917 - ---).
 \medskip \par \noindent
 $\bf 1912$ Birth of his son Ugo, a future prominent scientist in Physics. Student of Fermi and Heisenberg, he will be Professor Emeritus of the University of
 Chicago and Fermi award 1995.
 \medskip \par \noindent
 $\bf 1913$ Niels Bohr introduces a quantum model of the atom.
 \medskip \par \noindent
   $\bf 1914$ Elie Cartan writes an extended French translation of Fano's contribution to Enziklopaedie der Mathematische Wissenschaft on continuous groups.
  \medskip \par \noindent
  $\bf 1915$ Further papers by Fano, after those on cubic threefolds and \cite{56}, on three dimensional algebraic varieties with all plurigenera equal to zero.
  This theme, that is the theme of Fano threefolds, will accompany him for the rest of his life.
   \medskip \par \noindent
  $\bf 1914-18$ First World War. Engaged as a lieutenant, Fano is also office manager of the Piedmont Committee for industrial mobilization.\footnote{\tiny He was also editor of a book reporting the activities of the
  Committee during the war: \it L' opera del Comitato Regionale di Mobilitazione Industriale per il Piemonte (Settembre 1915 -  marzo 1919).\rm} For his service
  he was given the prestigious title of 'Ufficiale dell' Ordine della Corona d' Italia'. 
 \medskip \par \noindent
 $\bf 1917$ Birth of his son Roberto, a future prominent scientist in Computer Science. He is Professor Emeritus at  Massachussetts Institute of Technology.  Shannon 
  award 1976 for his work on Information theory and Electrical Engineering.
  \medskip \par \noindent
  $\bf 1922$ University of Torino: Fano gives the inaugural lecture for the beginning of the academic year. \medskip \par \noindent $\bf 1922$ Mussolini is prime minister after the fascist  coup 'Marcia su Roma'. \small
 \medskip \par \noindent $\bf 1923$ University College of Wales, Aberystwyth. Fano gives a series of invited lectures on the achievements of the Italian research in geometry and  two general lectures: 'Intuition in Mathematics' and 'All Geometry is theory of Relativity'\footnote{\tiny Parts of the corresponding manuscripts are contained in the Fondo Fano of the mathematical library of University of Torino. See also:  L. Giacardi  \it Testimonianze sulla Scuola italiana di geometria algebrica nei fondi manoscritti della Biblioteca  "Giuseppe Peano" di Torino \rm in 'Gli archivi della scienza. L' Universit\'a di Torino e altri casi italiani'  (S. Montaldo, P. Novaria editors)  Franco Angeli Milano (2011) 105-119.}
 \medskip \par \noindent
 $\bf 1925$ USSR: celebration of the 200 years of the Academy of Science. Gino Fano and Guido Fubini represent the University of Torino.
 \medskip \par \noindent
 $\bf 1926$ Kazan USSR: Fano gives a talk at the Conference for celebrating Lobachetvskij after 100 years of non euclidean geometries.
 \medskip \par \noindent
 $\bf 1928$ International Congress of Mathematicians, Bologna: two communications are presented by Fano. One is on Fano threefolds and their classification and the other one on
 plane birational contact transformations, see section 2.
 \medskip \par \noindent
 $\bf 1931$ All the professors of the Italian universities are required to pledge their loyalty to the fascist regime.\footnote{\tiny This happened with 12 exceptions.}
 \medskip \par  \noindent
 $\bf 1937$ All public clerks must be members of the National Fascist Party.
 \medskip \par \noindent
 $\bf 1928-38$ Several important papers are written by Fano on threefolds whose curvilinear sections are canonical curves, (Fano threefolds), and related topics.
 \medskip \par \noindent
 $\bf 1928-38$ Several contributions are offered by Fano to 'Istituto della Enciclopedia Italiana' and  to 'Enciclopedia delle Matematiche Elementari'. His book 'Geometria non euclidea (Introduzione geometrica alla teoria della relativit\'a) is published.
 \medskip \par \noindent
 $\bf 1938$ October 16: due to the racial laws against Jews, Fano is expelled from the ranks as a professor and from all Italian scientific institutions and academies.\footnote{ \tiny A list of them: R. Accademia delle Scienze di Torino, R. Accademia dei Lincei, R. Istituto Lombardo di Scienze e Lettere, R. Accademia Peloritana di Messina, Circolo Matematico di Palermo, Unione Matematica Italiana, Societ\'a Italiana per il Progresso delle Scienze, R. Accademia Virgiliana di Mantova.}
 \medskip \par \noindent
 $\bf 1939$ Fano and his wife finally decided to leave the country for Switzerland. His son Robert left few weeks later, reaching Bordeaux and then the U.S. in September 1939. The other son, Ugo Fano, came also there in the same year.
 \medskip \par \noindent
 $\bf 1939-45$ Refugee in Lausanne, Fano cooperates with the local University and the \'Ecole des Ing\`enieurs. He gives courses for the Italian university students
 of the refugees camps. Conferences at the 'Cercle math\'ematique de Lausanne'.
 \medskip \par \noindent
 $\bf 1940-45$ Papers in: Commentarii Mathematici Helvetici, Revista de la Universidad Nacional de Tucum\`an, the journals of Pontificia Academia Scientiarum.
 \medskip \par \noindent
 $\bf 1940-45$ Some topics: cubic threefolds and plane quintics. Cubic fourfolds with respect to rationality problem.\footnote{\tiny This is related to some results of Ugo Morin on the Pfaffian cubic fourfolds} Enriques surfaces (Reye congruences).  
 \medskip \par \noindent
 $\bf 1946$ Fano is back to Italy, restituted to his rights and status. He is retired and professor emeritus. He will spend most of his time in the U.S., visiting his sons.
 \medskip \par \noindent
  $\bf 1948$ Fano's last important memoir appears in Commentationes Pontificiae Acad. Sci.: it is on classification and rationality problems for Fano threefolds.   
  \medskip \par \noindent 
  $\bf 1950$ Torino, February: a conference in honor of Fano, including a series of five lectures by him, is organized at Seminario Matematico by his formers students.\footnote {\tiny For the list  of organizers see Rend. Sem. Mat. Univ. Polit. Torino 9 (1950) 5-7}
  \medskip \par \noindent
  $\bf 1952$ The last work of Fano was the preparation of the commemorative talk at  Accademia dei Lincei for his friend Guido Castelnuovo, dead on April 27 1952.
  \medskip \par \noindent
  $\bf 1952$ It will be impossible for him to give this talk. Fano passed the way before of the commemoration, in Verona on November 8, 1952.

   \newpage
   \section{Publications of Gino Fano 1890 - 1953}
\rm \small   In the next list of references we give the bibliographical coordinates of all papers written by Gino Fano we are aware of. The symbol * means that the paper is in the collection of scanned papers which follow this presentation.
   
\end{document}